\documentclass[12pt]{article}
\usepackage{amsfonts}
\usepackage{amsmath}
\usepackage{url}
\usepackage{stmaryrd}
\usepackage[all,cmtip]{xy}
\usepackage{tikz-cd}
\usetikzlibrary{shapes.geometric,arrows}
\begin{document}

\title{The Geometrization of Meaning}
\author{Michael Heller \\
Jagiellonian University \\
Copernicus Center for Interdisciplinary Studies\\
ul. Szczepa\'nska 1/5, 31-011 Cracow, Poland}

\date{\today}
\date{\today}
\maketitle

%\vspace{1.5cm}

\begin{abstract}
One of the greatest problems in philosophy is that of meaning. The turning point in thinking on meaning was Tarski's definition of truth, and the rapid development of logical semantics and model theory was a consequence of this achievement. Prhaps less well-known among classical logicians and philosophers is that it is category theory that provides the requisite mathematical tools to study the relationship between the syntax of formalized theories and their semantics. The aim of this article is to change this situation and make a preliminary philosophical analysis of the results obtained so far. They concern formalized algebraic theories with axioms in the form of equational laws, theories based on propositional logic and coherent Boolean logic, as well as decidable logic which is not necessarily Boolean. The syntactic-semantics relation for these theories takes the form of dualisms between the respective syntactic and semantic categories. These dualisms are given by the appropriate adjoint functors. In all the considered cases, the syntax--semantics dualism corresponds to the algebra-geometry dualism. We analyze the philosophical significance of these results. They allow us to look at the problem of meaning in a new light and formulate a criterion distinguishing formal theories from empirical theories. The disputes that were once fought over Tarski's definition of truth are transferred to a new context. The polemics as to whether Tarski actually succeeded in reducing the problem of meaning to purely syntactic terms has been superseded: between the syntactic and semantic categories there is not a relation of reducibility but rather that of interaction, and this relation is given by adjoint functors. We also touch upon other philosophical aspects of the categorical approach to the problem of meaning.
\end{abstract}

Keywords: Categorical Logic, Syntax-Semantics Duality, Meaning, Adjoint Functors.

\section{Introduction}
One of the most mysterious ``things'' in the world are the meanings of words. It is thanks to their meanings that words differ from mere sounds. Undoubtedly, many social processes and psychological mechanisms were involved in the creation of meanings, but meanings, when they arise, largely liberate themselves from the power of sociology and psychology: they create literature, science, and culture, which are to some extent autonomous. Are these phenomena not underpinned by some deeper logical or maybe even mathematical processes?

Modern physics was created when a few scientists at the turn of the 16th and 17th centuries managed to isolate from the thicket of phenomena occurring in the material world  a phenomenon that was simple enough to be subject to rational study, and yet complicated enough to carry essential information about the world. In this sense, the phenomenon of a stone freely falling in the gravitational field of the Earth created modern physics. What could be the equivalent of a freely falling stone in the analysis of the process of understanding?

The issue of meaning is situated somewhere at the interface between syntax and semantics, and is not just a cognitive science problem, but also underlies some key issues in natural science. For instance,
Paul Davies \cite[p. 167]{Davies} in the following way formulates ``the central problem of biology'': 
\begin{quote}
``what came first, complex organic chemistry or complex information patterns? Or did they somehow bootstrap each other into existence in lockstep? What is clear that chemistry alone falls short of explaining life. We must also account for the origin of organized information patterns. And not just information: we also need to know how \textit{logical operations} emerged from molecules, including digital information storage and mathematically coded instructions, implying as they do \textit{semantic} content. Semantic information is a higher-level concept that is simply meaningless at the level of molecules. Chemistry alone, however complex, can never produce the genetic code or contextual instructions.''\footnote{All of the italics are those of the author of this quote.}
\end{quote}
Let us notice that the ``organized information patterns'' and ``mathematically coded instructions'' are nothing else but syntactic issues, and to what they refer (their ``content'') is typically a subject-matter of semantic considerations. We can see how far this complex of issues is from the sought-after phenomenon of ``freely falling stone'', which would make this complex treatable. However, we do not think the problem is hopeless, and we can at least pave the way to its solution. For this purpose, the problem should be simplified as much as possible, that is, the problem of syntax should be reduced to that of a simple formal theory and its semantics should be recovered using strict methods of mathematical logic. This is a textbook issue, the novelty of our interpretation is that we base on the methods of mathematical category theory, which allow us to go much further than it was previously possible. Our approach differs from all the previous ones in that it not only examines problems arising in the interface between syntax and semantics, but also, using strict mathematical methods, shows how - on the one hand - syntax creates its semantics and - on the other hand - how syntax can be reconstructed from a given semantics. Moreover, it turns out that the relationship between syntax and semantics, in the case of relatively simple formalized theories, can be precisely described with the help of what category theorists call adjoint functors.

We have put the above quote from Paul Davies' book at the beginning of the present study not because we intend to address the problem raised by him, for we do not, but because it motivates well our intention to explore the relationship between syntax (information) and semantics,  and to do so with the help of innovative means provided by categorical logic. Although our considerations will be based on the primitive (as compared to the complexity with which we meet in biology) methods of mathematical logic and will be light yeaars away from the problem of the essence and origin of life, they may prove helpful in understanding how logic operations performed on the information are related to what the information is about.

In the first decades of the twentieth century, the syntactic approach to science predominated and semantic concepts were treated with suspicion. This was largely the result of the great authority enjoyed by Hilbert's program in mathematics, and then of the strong influence of the philosophy of logical empiricism. The situation changed radically under the influence of Tarski's work on the definition of truth \cite{Tarski}. Under the influence of Tarski, there was a lush development of semantics and model theory. Since then, both logic and metamathematics (which abundantly benefits from the achievements of the former) have made great strides. Thanks to the mathematical methods of category theory, the relations between syntax and semantics have become the area of geometric constructions and proving of theorems. The results achieved so far are so significant that they call for a philosophical analysis. This is the purpose of the present study.  However, since these results have been achieved by using advanced mathematical methods and they presuppose  expert knowledge in category theory and categorical logic, they cannot be presented in a popular way. Although we try to present them in as simple a manner as possible, reading our text would require considerable effort and mental discipline from a less prepared reader.

We begin, in section 2, by recalling some basic concepts from the categorical logic that will be necessary for our further analysis. In the next five sections, we review four approaches to the topic of interest, in order of increasing generality and advancement of mathematical tools. And so in section 3, we present a result that can already be considered classic: if a formalized theory assumes only propositional logic, then there is a duality between its syntax and semantics, which is expressed by the appropriately reformulated Stone theorem. In setion 4, we give an interpretation of this result. We can go further (section 5) if we assume that the theory under consideration is an algebraic theory (such as, for example, group theory); such a theory presupposes first order logic and its axioms assume the form of equational laws. It was this case that initiated a categorical version of the whole issue \cite{Lawvere}. In this case, the relationship between syntax and semantics takes the form of a dualism between the category of algebraic theories and the category of their models. In section 6, we present what is known as the Makkai duality \cite{Makkai1,Makkai2,Makkai3}. It is the duality between the theories presupposing coherent Boolean logic and the space of models of such theories organized into the form of the so-called ultragrupoid. An alternative version of this case was proposed by Awodey and Forssell \cite{AwodeyForssell,Forssell}. In their approach (section 7), a  coherent category need not be Boolean as long as it is decidable (admitting inequality). On the semantic side, instead of the ultragroupoid, there is a topological groupoid with the so-called logical topology. In the last two sections, we put the whole thing into a broader perspective. In section 8, we give brief comments on the reported results with special emphasis on how they change our understanding of the ``meaning of meaning''. In section 9, we refer to the old discussion of whether semantics can be reduced to syntax, and signal the possibility of applying our results to non-formal situations.

\section{Syntactic Category}
We understand a theory, in broad terms, as a formal language adapted to axiomatize a certain domain of knowledge. The concept is broad enough to cover notion of theory as it is functioning in physics (although in physical practice theories are seldom fully formalized), but sufficiently formal to be itself subject to a strict mathematical analysis. However, not to simplify too much, we will consider many-sorted theories, i.e., theories the signature $\Sigma $ of which gives a finite list of types (sorts), and all expressions are qualified as being of a given type. \textit{Univers du discours} of a many-sorted theory is not homogeneous but contains entities of various kinds (sorts or types) which are subject to different logical manipulations. For instance, in elementary Euclidean geometry one considers points, lines and planes (as different types) which, even in non-formal presentations, are distinguished by different kinds of symbols. Let us now be more technical.

A signature $\Sigma $ of a formal language ${\cal L}$ contains a list of symbols: predicate or relation symbols (such as $ \in , \leq )$ , function symbols (such as $+, \times , 0, 1$), constant symbols and variables. $\Sigma $ must contain at least one type symbol.
\begin{itemize}
	\item Each predicate symbol $p \in \Sigma $ has an arity $\sigma_1 \times \ldots \times \sigma_n$ with $\sigma_1, \ldots , \sigma_n \in \Sigma $ being type symbols (they need not be different).
	\item Each function symbol $f \in \Sigma $ has an arity $\sigma_1 \times \ldots \times \sigma_n \to \sigma $ with $\sigma_1, \ldots , \sigma_n, \sigma \in \Sigma $.
	\item Each constant symbol $c \in \Sigma $ and each variable $x, y, \dots \in \Sigma $ has its type assigned.
\end{itemize}

Let us consider a language ${\cal L}$ specified by its signature $\Sigma $. By a theory $T$ we understand the collection of assertions, which could be derived -- by using deduction rules of a given logic -- from a chosen set of other assertions, called axioms of this theory (for precise definition see \cite{BarrettHalvorson}, \cite[p. 527-530]{MacLaneMoerdijk}). Different systems of axioms may lead to the same theory.

So far everything is well known from the logical studies of formal languuages. Now, we want to employ powerful tools provided by category theory. The starting point is that any formal theory $T$, as it was characterized above, can be organized to form a category. It is done in the following way.

Let $\Phi $ be a collection of all well-defined formulae of $T$. The pair $(\Sigma , \Phi )$ is called a context. $\Sigma $ contains a collection of rules assigning a type to each term of the theory $T$. In this formalization, the pair $(\Sigma , \Phi )$ indeed plays the role of a context for well defined formulae of $\Phi $.

The syntactic category Syn($T$) of a formal theory $T$, also called a category of contexts,  is defined as follows. 

%The following paragraph should be expllained.
Objects of Syn($T$) are contexts $(\Sigma , \Phi )$; less formally, an object of Syn($T$) is a formula $\varphi(\bar{x}) \in \Phi $, where $\bar{x} = (x_1, \ldots x_n)$, up to renaming free variables.

Let now $\varphi(\bar{x})$, $\bar{x} = (x_1, \ldots x_n)$, and $\psi(\bar{y})$, $\bar{y} = (y_1, \ldots y_m)$ be two objects of Syn($T$). We assume (without the loss of generality) that the collections $\varphi(\bar{x})$ and $\psi(\bar{y})$ are different. We say that a functional relation $\chi(\bar{x}, \bar{y})$ from $\varphi $ to $\psi $ is defined if from the theory $T$ it follows that, for each $\bar{x}$ such that $\varphi(\bar{x})$, there is a unique $\bar{y}$ such that $\psi(\bar{y})$, and $\chi(\bar{x}, \bar{y})$.

The morphisms $(\Sigma , \Phi ) \to (\Pi , \Psi )$ of Syn($T$) are interpretations of variables. To be more precise, a morphism from  $\varphi \in \Phi $ to $\psi \in \Psi $ is an equivalence class $[\chi ]$ of functional relations from $\varphi $ to $\psi $, where $\chi_1 $ and $\chi_2 $ are equivalent if in $T$ one has
\[
\chi_1 (\bar{x}, \bar{y}) \vdash \chi_2 (\bar{x}, \bar{y}),
\]
and automatically,
\[
\chi_2 (\bar{x}, \bar{y}) \vdash \chi_1 (\bar{x}, \bar{y}),
\]
for details see, e.g. \cite{Oosten}.

This is a highly condensed description of what these constructions really look like in  textbooks of logic. However, one loophole has to be filled in for this description to make sense at all. We have to determine the logic on which we make all our inferences. Moreover, the laws of this logic must be sanctioned by axioms that appear in our constructions. In the following sections, we will look at how logic organizes the interaction between syntax and semantics.

\section{Stone Duality and Semantics}
To see how our scheme works in practice, let us assume that the logic in question is ordinary propositional logic. Consequently, the theory $T$ we are going to consider is a theory expressed in terms of this logic. It is one-sort theory, and its signature is $\Sigma = \{p_0, p_1, \ldots \}$, where $p_0, p_1, \ldots $ stand for propositions (sentences).\footnote{To refresh your knowledge of propositional logic consult \cite[chapter 6]{Goldblatt}.} Let us see what the syntactic category looks like in this case.

In the case  of propositional logic, a formula in context reduces to just a sentence and, consequently, sentences are objects of Syn($T$). If $\varphi $ and $\psi $ are two such objects, there is exactly one (up to equivalence) functional relation from $\varphi $ to $\psi $, namely $\varphi \vdash \psi $. And this is precisely the morphism from $\varphi $ to $\psi $ (see \cite{HalvTsem}).

If we remember that the Lindenbaum-Tarski algebra for a theory $T$ consists of the equivalence classes of sentences of $T$ that are provably equivalent in $T$,\footnote{I.e., two sentences are provably equivalent if each of them implies the other.} we then easily conclude that the syntactic category Syn($T$) for a propositional theory $T$ is the same thing as the Lindebaum-Tarski algebra of $T$.

Logical operations in any Lindenbaum-Tarski algebra respect the above mentioned equivalence relation, and any such algebra has naturally built-in order relation in the following way
\[
[\varphi ] \leq [\psi ] \Leftrightarrow T \vdash \varphi \to \psi .
\]
This makes it a Boolean algebra\footnote{On Boolean algebras. \cite[pp. 133-135]{Goldblatt}}. And conversely, with any Boolean algebra $B$ we can associate a propositional theory $T_B$. We start by constructing a language $\cal{L}_B$ for this theory: with each element $b \in B$ we associate a propositional constant (representing a particular sentence). And then we construct the theory $T_B$ essentially by translating $b \leq c$, $b, c \in B$ into the implication of the corresponding propositional constants, and analogously for other logical operations (for details see \cite[section 1.2.1]{Forssell}).

Two propositional theories are algebraically equivalent if their Lindenbaum-Tarski algebras are isomorphic. Two such theories can be regarded as expressing the same content in two different languages.

Let us proceed further. A subset $F$ of a Boolean algerbra $B$ is said to be a filter if
\begin{itemize}
\item $1 \in F$ and $0 \not\in F$,
\item if $a \leq b$ and $a \in F$ then $b \in F$,
\item if $a \in F$ and $b \in F$ then $a \wedge b \in F$
\end{itemize}
Filters are partially ordered by inclusion. A filter on $B$ that cannot be enlarged to a bigger filter on $B$ is called an ultrafilter. It can easily be seen that a filter $F$ is an ultrafilter if, for every $b \in B $, either $b \in F$ or $\neg b \in F$.

Let us now define the set $D(b)$, $b \in B$, of ultrafilters  by $U \in D(b) $ iff $b \in U$. The sets $D(b)$, for all $b \in B$, form a basis for a topology on the space of all ultrafilters in $B$. Thus its points are ultrafilters of $B$.  This space is called the Stone space, denoted Stone($B$).\footnote{For details see \cite[pp. 515-516]{MacLaneMoerdijk}. The Stone space can also be defined as the space of prime ideals in $B$, since each ultrafilter $U$ determines a prime ideal ${\cal P} = \{b| \neg b \in U\}$, and \textit{vice versa}.}

Boolean algebras form a category, denoted {\bf BA}. Its objects are, of course, Boolean algebras and its morphisms Boolean homomorphisms, i.e. functors between Boolean algebras that preserve the structure of Boolean algebras (see \cite[p. 35]{Awodeybook}). The Boolean algebra ${\bf 2} = \{0, 1\}$ is the initial object of this category, and can be regarded as a ``set of truth values''. If so, the morphism  $p: B \to {\bf 2}$ is but a model of a Boolean algebra $B$.

On the other hand, any morphism from any object of ${\bf BA}$ to the initial object ${\bf 2}$, $p: B \to {\bf 2}$ gives an ultrafilter $U_p = p^{-1}(1)$, and \textit{vice versa} given an ultrafilter $U \subset B$, one obtains $p_U: B \to {\bf 2}$ by taking $p_U(b) = 1, \; p \in B$. Consequently, a Stone space (its points being ultrafilters) can be regarded as a space of two-valued models of a Boolean algebra seen as a propositional theory modulo the ``provability equivalence''. Therefore, we obtain the space of models Mod(B)  of the theory $T$ that corresponds to the Boolean algebra $B$ as the set of morphisms from the Boolean algebra $B$ to the Boolean algebra ${\bf 2}$, i.e.
\[
{\rm Mod}(B) = \mathrm{Hom}_{\bf BA}(B, {\bf 2}).
\]

Can we retrieve from this space of models the original Boolean algebra $B$, and consequently the corresponding theory $T$? The answer is ``yes'' provided we equip Mod($B$) with a suitable topology, the so-called ``logical'' topology. Its basic open sets are given in terms of sentences $\varphi $ of the language ${\cal L}$ of $T$
\[
V_{[\varphi ]} = \{Y \in {\rm Mod}(B)| Y \models \varphi\}.
\]
Then the Boolean algebra $B$ is recovered by taking morphisms in the category ${\bf Stone}$ from Mod($B$) to ${\bf 2}$
\[
B = \mathrm{Hom}_{\bf Stone} (M_B, {\bf 2}).
\]

As we can see, Mod($B$), regarded as a mere collection of models of a theory $T$, contains less information than $T$ itself. The lacking information has to be supplemented by providing a suitable topology on Mod($B$). This means that the theory $T$ contains information about how its models interact topologically with each other, while the mere collection of models Mod($B$) does not contain such information \cite[pp. 411-413]{HalvTsem}.

In fact, we have here a functorial dependence given by the following pair of functors

\begin{center}
\begin{tikzcd}
\mathbf{BA}^{\rm op} \arrow[rrrrr, "{{\rm Hom}_{\bf BA}(-,{\bf 2})}", bend right] &  &  &  &  & \bf Top \arrow[lllll, "{{\rm Hom}_{\bf Top}(-, {\bf 2})}"', bend right]
\end{tikzcd}

\small{Diag. 1}
\end{center}
which form the adjunction. However, if we restrict the category ${\bf Top}$ to its subcategory ${\bf Stone}$ of Stone spaces (and continuous functions), the adjunction changes into the equivalence of categories (for details and proofs see \cite[section 1.2.1]{Forssell}). If this is the case, the above diagram assumes the form

\begin{center}
\begin{tikzcd}
\mathbf{BA}^{\rm op} \arrow[rrrrr, "{{\rm Hom}_{\bf BA}(-,{\bf 2})}", bend right] &  &  &  &  & \bf Stone \arrow[lllll, "{{\rm Hom}_{\bf Stone}(-, {\bf 2})}"', bend right]
\end{tikzcd}

\small{Diag. 2}
\end{center}
and the adjunction changes into the equivalence of categories.

\section{Interpretation}
Diag. 2 summarizes what is known as the Stone duality theorem asserting that there exists a (contravariant) equivalence of the category ${\bf BA}$ of Boolean algebras and the category ${\bf Stone }$ of Stone spaces. On the left hand side of Diag. 2 we have the category ${\bf BA}$ of Boolean algebras, which in fact represents the content of a certain propositional theory $T$; on its right hand side we have the category ${\bf Stone}$ of Stone spaces which is but a space of models of the theory $T$ modulo its linguistic expression. The left hand side is clearly algebraic and syntactic, whereas the right hand side is geometric (as a space) and semantic (as a space of models). The fact that this space is equipped with a ``logical'' topology guarantees that provable sentences of the theory $T$ are  true in all models \cite{Forssell}.

The interaction between syntax and semantics is given by the functors ${\rm Hom}_{\bf BA}(-,{\bf 2})$ and ${\rm Hom}_{{\bf Top}({\bf Stone})}(-,{\bf 2})$. In the case of ${\bf Top}$ (Diag. 1) we have the adjunction, and in the case of ${\bf Stone}$ (Diag. 2) the equivalence of categories. As it is well known, it is the equivalence of categories that determines the ``mathematical identity'' of a category. A functor $F: {\cal C} \to {\cal D}$ defines an equivalence between the categories ${\cal C}$ and ${\cal D}$, if there exists a functor $G: {\cal D} \to {\cal C}$ such that there are natural isomorphisms ${\bf 1}_{\cal C} \equiv G \circ F$ and ${\bf 1}_{\cal D} \equiv F \circ G$. The requirement that the images of the compositions $F \circ G$ and $G \circ F$ should naturally transform into the respective identities, leaves quite a lot of room for both categories to differ from each other. For instance, the ``size'' of these two categories can be very different (\cite[pp. 161-162]{Marquis}). We thus see that even in such a simple case as is exemplified by propositional logic, the interaction between syntax and semantics is not just of one-to-one type, but requires the subtle machinery of category theory to be fully articulated.

To grasp the intuition behind this approach, let us once more consider Diag. 2. It doubtlessly has a geometric aspect (as a figure on the Euclidean plane); however, it was produced by the following program

\begin{verbatim}
\begin{tikzcd}
\mathbf{BA}^{\rm op} \arrow[rrrrr, "{{\rm Hom}_{\bf BA}
(-,{\bf 2})}", bend right]  &  &  &  &  & 
\bf Stone \arrow[lllll, "{{\rm Hom}_{\bf Stone}
(-, {\bf 2})}"', bend right]
\end{tikzcd}
\end{verbatim}

\noindent which is clearly linguistic in nature, it consists of a string of symbols subject to certain syntactic rules. It can be regarded as a ``theory'' to produce graphics. In this particular case, this ``theory'' \textit{is about} graphic given by Diagram 2. We thus have a nice example of the duality between syntax and semantics -- algebra and geometry. Let us notice, however, that Stone's theorem speaks of a duality between two categories, while in our example we were dealing with a single theory and some  semantics of it. In order to present our example more precisely, we would first have to define the appropriate category of programs and the appropriate category of diagrams, and then prove the existence of a certain duality between them (and, if necessary, to enrich the required logic). Nevertheless, our example, in its present form, nicely shows that the syntactic-semantic interpretation of Stone's theorem is in line with our intuitions. Falling stone (see Introduction) created modern physics, will Stone's theorem play the same role in the geometric theory of meaning?

\section{Algebraic Theories}
Propositional theories are a good first step in analyzing the interplay between syntax and semantics. The fact that this is a simple (yet not trivial) case made it possible to isolate the issue of interest from among many entanglements and complications. However, the next step has to be taken. It is provided by the so-called algebraic theories, known also as equational theories. The idea and its seminal implementation was proposed by Lawvere in his influential PhD thesis \cite{Lawvere}.

Algebraic theories are first order theories, the signatures of which have function symbols but no relation symbols (except that of equality), and all their axioms are equational laws. Examples are: groups, rings, moduli.\footnote{One important thing is that operations should be defined everywhere. This excludes, for example, fields since the inverse of zero is not defined \cite[section 2.1]{AwodeyBauer}}

We now construct the syntactic category Syn($A$) of an algebraic theory $A$ in the way described in section 2 adapted to peculiarities of algebraic theories \cite[section 2.1.3]{AwodeyBauer}. Syn($A$) is a small category with finite products. Its objects  form a sequence $A^0, A^1, A^2, \ldots $ such that $A^m \times A^n = A^{m+n}$, $m, n \in \mathbb{N}$, and its morphisms are finite lists of equivalence classes of terms, provable in $A$; for instance
\[
[f_1, \ldots , f_m]: T_n \to T_m
\]
where each $f_i$ has arity $1$. Moreover, $1 = A^0$ is the terminal object, and each  object is a product of finitely many copies of $A^1$. 

And \textit{vice versa}, having a category Syn($A$), we can recover the corresponding algebraic theory $A$ by regarding morphisms of Syn($A$) as function symbols of $A$, and interpreting equalities between morphisms of Syn($A$) as axioms of $A$. Two algebraic theories are categorically equivalent if their syntactic categories are equivalent.

The next move is to consider the category ${\cal FP}$ of finite product categories with finite product preserving functors as morphisms. The category ${\cal E}$ of algebraic categories is a subcategory of ${\cal FP}$. The category {\bf Sets} is an object of ${\cal FP}$. An important role in our construction is played by the category ${\cal G}$ of categories with all limits and colimits, and functors preserving limits, filtered colimits\footnote{A category $C$, in which every finite diagram has a cocone, is called a finitely filtered category. A diagram $F: D \to C$, where $D$ is a finitely filtered category, is called a filtered diagram. A colimit of a filtered diagram is said to be a filtered colimit.} and regular epimorphisms as morphisms. The category ${\bf Sets }$ is also an object of ${\cal G}$. The functor ${\rm Hom}_{\cal FP}(-, {\bf Sets})$ is, in fact, a contravariant functor from ${\cal FP}$ to ${\cal G }$, and it has a right adjoint functor ${\rm Hom}_{\cal G}(-, {\bf Sets})$ from ${\cal G}$ to ${\cal FP}$ (for details see \cite[section 1.2.2]{Forssell}) .

A set-valued model of an algebraic theory $A$ is a finite product preserving functor from Syn($A$) to ${\bf Sets }$. And we have the category Mod($A$) of set-valued models of $A$
\[
{\rm Mod}(A) \cong {\rm Hom}_{\cal FP}({\rm Syn}(A), {\bf Sets })
\]
with morphisms as natural transformations of such functors.

If we now regard algebraic theories as objects of the category ${\cal E}$, we can restrict the functor 
\[
{\rm Hom}_{\cal FP}(-, {\bf Sets}): {\cal FP}^{op} \to {\cal G}
\]
to ${\cal E}$, and consider the image ${\cal MOD}$ of the above functor restricted to ${\cal E}^{op}$. ${\cal MOD}$ is clearly the category (with morphisms as  natural transformations of corresponding functors suitably restricted) of models of algebraic theories as objects of the category ${\cal E}$.

And now we have all essential concepts ready to formulate the final result: \textit{There is an adjunction between the category of algebraic theories and the category of their models as shown in the following diagram}

\begin{center}
\begin{tikzcd}
\mathbf{\cal E}^{\rm op} \arrow[rrrrr, "{{\rm Hom}_{\cal E}(-,{\bf Sets})}", bend right] &  &  &  &  & {\cal MOD} \arrow[lllll, "{{\rm Hom}_{\cal MOD}(-, {\bf Sets})}"', bend right]
\end{tikzcd}

\small{Diag. 3}
\end{center}

\noindent For proof see \cite[Theorem 1.2.2.1]{Forssell}. The interpretation of this result is clear: algebraic theories generate the entities they speak of (their semantics), but these entities are not something external to these theories, but are adjointly coupled to them. 

\section{Ultragroupoid Semantics}
It is understandable that as we enrich the theory under consideration, the construction of semantics becomes more and more complicated. The degree of complication grows so quickly that, already when going to the next step, we will have to be content with an even more simplified way of presentation. We will limit ourselves only to sketching  constructions, merely giving some hints in the footnotes towards the definition of concepts.

Our next step consists of replacing propositional theories of section  3 by  theories  employing the so-called coherent logic, and the Lindenbaum-Tarski algebra of a given propositional theory $T$  by the Boolean coherent category. Coherent logic is a fragment of (finitary) first-order logic admitting only the connectives (and, or, true and false) and the existential quantifier. We additionally allow equality. Syntactic category of such a theory is a coherent category. The latter is a regular category\footnote{A regular category is a category the internal logic of which is a regular logic. Logical operations of this logic consist only of truth, conjunction, and existential quantifier.} the subobject posets Sub($X$) of which have finite unions, and these unions are preserved by the functors $f^*: {\rm Sub}(Y) \to {\rm Sub}(X)$. A Boolean coherent category is a coherent category in which every subobject has a complement.\footnote{This means that for any monomorphism $A \hookrightarrow X$ there is a monomorphism $B \hookrightarrow X$ such that $A \cap B$ is the initial object and $A \cup B=X$. The name ``Boolean'' comes from the fact that the subobject lattice Sub($X$) of any object $X$ of a Boolean coherent category is a Boolean algebra.} We can thus simply identify the category of first order theories with the category of the Boolean coherent categories (with coherent functors as morphisms.\footnote{A functor between coherent categories is a coherent functor if it is a regular functor and preserves finite unions.}) And, if a Boolean coherent category has some additional properties,\footnote{It should have: stable, disjoint finite coproducts and coequalizers of equivalence relations.} it is a Boolean pretopos.\footnote{If we remove the postulate of the existence of power objects from the definition of topos, we get pretopos. A pretopos, which is also a Boolean category, is a Boolean pretopos.}

In an analogous way, as morphisms $B \to {\bf 2}$ are models of propositional theories, coherent functors from a Boolean coherent category to the category of sets, or from Boolean pretopoi to the category of sets, are set-valued models of first-order theories.

Makkai was able to show that there exists the adjunction between the category ${\cal BP}$ of Boolean pretopoi and the category ${\cal UG}$ of ultragroupoids \cite{Makkai1,Makkai2,Makkai3}. The most condensed groupoid definition is that it is a category in which every morphism is invertible, and ultragroupoid is a groupoid with a kind of ultraproduct structure.\footnote{If we consider a family of structures of the same signature, we can form a new structure of the same signature with the help of ultrafilters. This new structure is called an ultraproduct (this is not a definition!).} This adjunction has the form

\begin{center}
\begin{tikzcd}
\mathbf{\cal BP}^{\rm op} \arrow[rrrrr, "{{\rm Hom}_{\cal BP}(-,{\bf Sets})}", bend right] &  &  &  &  & {\cal UG} \arrow[lllll, "{{\rm Hom}_{\cal UG}(-, {\bf Sets})}"', bend right]
\end{tikzcd}

\small{Diag. 4}
\end{center}

As we can see in this diagram, the functor ${\rm Hom}_{\cal BP}(-, {\bf Sets})$ sends a Boolean pretopos ${\cal B}$ to its ultragroupoid ${\rm Hom}_{\cal BP}({\cal B}, {\bf Sets})$ and, conversely, the functor ${\rm Hom}_{\cal UG}(-, {\bf Sets})$ sends any ultragroupoid $K$ to the corresponding Boolean pretopos ${\rm Hom}_{\cal UG}(K, {\bf Sets})$. It can be shown that a pretopos (representing a first order theory) can be recovered from its groupoid of models.

\section{Continuous Groupoid Semantics}
An alternative approach to Makkai's work was proposed by Awodey and Forssell \cite{AwodeyForssell,Forssell}. Their approach is more geomerical in character than that of Makkai, and in the case of propositional theories it specializes to the Stone duality. On the syntactical side it also considers Boolean coherent categories, but  can include non-Boolean coherent categories as well under the condition that they are decidable. Coherent decidable categories\footnote{An object $X$ of a coherent category ${\cal C}$ is called decidable if its  diagonal morphism $\Delta X: X \to X \times X$ is complemented. This implies that, for an object X,  in the internal logic of the category ${\cal C}$, it is true that  either $x=y$ or $ x \neq y$ for any $x,y \in X$.} represent coherent theories with the inequality $(\neq )$ predicate for each type.

On semantical side, in this approach, we have a space of models of the theory $T$ just as was the case with the propositional theories, except that now this space is more complicated. One then considers the category of models Mod($T$) of a given Boolean coherent (or non-Boolean decidable) theory $T$, where a model of $T$ is a functor from the categorical (algebraic) representation ${\cal C}_T$ of $T$ into ${\bf Sets }$, i.e.
\[
{\rm Mod}(T) \simeq {\rm Hom}({\cal C}_T, {\bf Sets})
\]
with the category ${\bf Sets}$  playing the role of a dualizing object, just like ${\bf 2}$ in the Stone duality.

If we eliminate from Mod($T$) all non-isomorphism arrows, we obtain a groupoid, i.e. a category in which every morphism is an isomorphism. When studying the Stone duality we had to equip Mod($T$) with the ``logical'' topology, defined in terms of sentences of the theory $T$, now we must generalize this topology with the help of much richer structure, namely we must consider a Grothendieck topos of equivariant sheaves on a topological groupoid of models.\footnote{Let us consider a space $X$ with action $\rho : G\times X \to X$, where $G$ is a group, and projection $p: G \times X \to → X$. These data lead to an action groupoid in the category of spaces. An equivariant sheaf over a space $X$ is a sheaf $x$ over $X$ with an isomorphism $\theta :p^*x \to ^*x$ of sheaves over $G \times X$ (with the cocycle condition on $G \times G \times X$ \cite{nLabEquiv}).} An explanation is necessary here.

Every Grothendieck topos is a topos of sheaves, and every sheaf can be thought of as a generalized topology. A family of open sets ${\cal O}(X)$ defining a topology on a space $X$ can be regarded as a category with open sets from ${\cal O}(X)$ as objects and inclusions between open sets as morphisms. Let us now consider a functor category ${\bf Sets}^{{\cal O}(X)^{op}}$. Objects of this category are called presheaves and, with natural transformations between them as morphisms, form a category. The full subcategory of this category, together with the condition guaranteeing that functions defined locally could be ``glued together'' in a global sense, is called the category of sheaves, denoted  by Sh($X$).\footnote{In other words, we apply here the operation of ``sheafification'' that turns presheaves into sheaves \cite{nLabSheafif}.}  In this way, any Grothendieck topos can be regarded as a generalized topological space in that it allows to ``study things locally that are transformed continuously into one another'' \cite[p. 163, 249]{Marquis}.\footnote{In category theory, a Grothendieck topology is a structure on a category ${\cal C}$ that makes the objects of ${\cal C}$  act like the open sets of a topological space. A category together with a choice of Grothendieck topology is called a site. Every Grothendieck topos arises as sheaves on a site.} Now, we topologize the groupoid of models $G = {\rm Mod}(T)$ with the help of Sh($G$).

If we apply the above procedure to Mod($T$), it becomes a topological groupoid of models (with isomorphisms between them as morphisms), equipped with the logical topology.

Now, it can be proved \cite{AwodeyForssell,Forssell} that this structure on Mod($T$) is enough to recover from it (up to a form of equivalence) the original theory $T$. This implies that, if  categories Mod$({\cal C}_T, {\bf Sets})$ and Mod$({\cal C}_{T'}, {\bf Sets})$ are equivalent, then the syntactic categories ${\cal C}_{T}$ and ${\cal C}_{T'}$ are also equivalent, meaning in turn that the theories $T$ and $T'$ are syntactically equivalent, in a suitable sense.

We thus have, on the syntactic side, decidable (not necessarily Boolean) category; let us denote it by {\bf dCoh}, and on semantic side a topological groupoid category. Is it possible to determine the kind of adjunction  between them? The answer is in the affirmative as long as we impose strict conditions on these structures. Since the formulation of these conditions would break the framework of this article, we will limit ourselves to the final result only (for details see \cite{AwodeyForssell}). 

First, we must narrow our considerations to the category {\bf wcGpd} of weakly coherent groupoids (and compact morphisms) \cite[Deinition 2.1.16]{AwodeyForssell} and we define two functors: the contravariant semantic functor
\[
{\rm Mod}: {\bf dCoh}^{op} \to {\bf Gpd},
\]
which factors through {\bf wcGpd}, and the semantic functor
\[
{\rm Form}: {\bf wcGpd} \to {\bf dCoh}^{op}.
\]
Then we show that they are adjoint functors \cite[Theorem 2.5.3]{AwodeyForssell}. This leads to the diagram

\begin{center}
\begin{tikzcd}
{\bf dCoh}^{op} \arrow[rrr, "{\rm Mod}"', bend left] &  &  & {\bf wcGpd} \arrow[lll, "{\rm Form}"', bend left]
\end{tikzcd}

\vspace{0.30cm}
\small{Diag. 5}
\end{center}
It can also be demonstrated that if ${\cal D} \in {\bf dCoh}$ is a pretopos, this adjunction changes into the equivalence of categories.

In this approach, the duality between decidable coherent theories and topological groupoids of models cleary exhibit an algebra--geometry duality aspect: syntax takes the algebraic attire, while semantics assumes the geometric shape.

\section{Analysis}
Let us start with a simple example. The utterance ``Someone is writing a letter'' is not a sentence in the logical sense since we cannot decide whether this utterance is true or false as long as we do not know who that someone is. As soon as we find out that this is John, we make an \textit{interpretation} of this statement, and if it turns out that this interpretation leads to a true sentence, we say that we have obtained a \textit{model} of this sentence; in other words, that ``John'' satisfies this sentence. We see that the concept of model, as a morphism from a Boolean algebra to category ${\bf 2}$, that occurs in our presentation of the Stone dualism, is simply an algebraic-categorical stylization of the above example. 

At the heart of all this is Tarski's observation that what a language refers to can be defined without going beyond the strict limits of formalism. All works reviewed in the preceding sections remain within this ``Tarskian paradigm''. An important novelty that these works introduce to this paradigm consists in showing the dualism (in a mathematical sense) between syntax and semantics, and that it is a kind of dualism between algebra and geometry with algebra playing the role of a language (syntax) and geometry of what this language is about (semantics). The latter dualism has been known, albeit not in formalized form, at least since Descartes' analytic geometry, but the discovery that it is related to the syntax -- semantics interaction is something relatively new.

As we have seen in section 5, this dualism works for those mathematical theories which can be formalized as algebraic (equational) theories. The logic presupposed by these theories is the ordinary first-order logic augmented with equality as the only predicate symbol. The model theory for this logic was earlier worked out by Birkhoff, Malcev and others (for the history see \cite{Elgueta}), and developed into the universal algebra which studies algebraic structures \textit{in abstracto}. For instance, it takes the class of groups as an object of its studies rather than concrete groups. This research thread was taken up by Lawvere in his seminal doctoral thesis \cite{Lawvere}. The essence of his idea consists in developing the concept of an algebraic theory (as outlined above) and elaborating mathematical tools with the help of which this concept could become the subject of a rigorous mathematical investigation. The problem comes down to finding a form of algebraic theory that would express its true content independent of any linguistic particulars (such as signature and specific form of axioms). This is similar to the problem of relativity theory in physics: how to find a form of a physical law that would be independent of any reference frame  with the proviso that in the present context the ``frame of reference'' is to be understood so that what does not depend on it is the \textit{meaning} of an algebraic theory. Lawvere showed that such an "invariant" form of algebraic theory is its categorical formulation. Another credit to Lawvere is that he was able to see that the ``paradigmatic method'' to study this problem is to use adjoint functors \cite[p. 194]{Marquis}. The adjunction between the category of algebraic categories and the category of their models (diagram 3) displays, in fact, a duality between ``syntactically invariant'' algebraic theories and what these theories speak about.

In light of the above, it is no longer true that ``mathematics may be defined as the subject in which we never know what we are talking about, nor whether what we are saying is true'' (Bertrand Russell). Mathematics is about its own ``semantics'' that is an essential part of itself. Strictly speaking, this is true only with regard to algebraic theories (since it has only been proven with regard to them), but we can presume that analogous theorems could  be proved also for some other formalized mathematical theories.

Both Makkai's and Awodey--Forssell's approaches have similar philosophical significance. They extend the syntax-semantics duality from propositional theories to first-order theories, with the Awodey-Forssell approach going a step further by allowing the inclusion of nonecessarily Boolean theories into the scheme. In both approaches, the space of models has the groupoid structure: in Makkai's approach it is an ultragroupoid, whereas in the case of Awodey and Forssell it is a topological groupoid with the topology given by a Grothendieck topos.

In all four cases considered in Sections 3 - 5, we have observed a strong interaction of the syntactic structure of a given theory (in its ``invariant'' version) with the space of its models. But this interaction comes fully to the fore (via adjoint functors) only when the  space of models is properly structured. And what might come as a surprise is that a simple collection of all models of a given theory contains less information than the theory itself. The theory potentially contains all its models, however not as a loose set, but as an organized totality, that is, it also takes into account all possible interactions between the models. In order for the space of models to include this information, it must be enriched with an appropriate geometric structure; in our case, either as an ultragroupoid, or as a topological groupoid with an appropriately rich topology.

All theories considered in our study are formal theories. According to standard philosophical doctrine, such theories consist of tautologies, and tautologies ``speak of nothing''. This doctrine goes back to Wittgenstein's \textit{Tractatus Logico-Philosophicus}. From Theses 4.461 and 4.462 we learn that the meaning of the proposition is determined by its verifiability. However, tautologies and contradictions are beyond verification because the former are always true and the latter are always false. Consequently, they have no meaning. In other words,  tautologies are propositions that are true (in the Tarskian sense) in any of their domains, i.e. of which every domain is a model. And contradictions are propositions that are false in any of their domains, i.e. that have no models. 

If, instead of propositions, we talk about those formal theories (i.e. sets of propositions) that were analyzed in the previous sections, then we are able to define the geometric structure of the space of all models of a given theory (Stone space, groupoid). Of course, for contradictions the space of models is empty. We can go a step further and distinguish a certain subspace of the space of models of a given theory $T$ as models that are consistent with some kind of experiment. The models belonging to this subspace deserve to be called the physical or empirical models of the theory $T$. This suggests that the space of models determines the meaning of a proposition or theory. From this perspective, tautologies are not meaningless, but -- on the contrary -- full of meaning. In a sense, however, this meaning is trivial since, for all practical purposes, there is no essential difference between the Wittgensteinian view that tautologies are devoid of meaning and the above view that tautologies are overflown with meaning, just as, from the point of view of pouring wine into a glass, there is no difference between a bottomless glass and a glass filled to the brim. In both cases, nothing can be changed in the state of affairs. In this view, the meaning of empirical theories is nontrivial; it is determined by the ``size'' of the distinguished subspace of the space of models. In a particular case, if we define an appropriate measure in the space of models, we could mathematically determine the ``degree of empiricality'' of a given theory, and we could speculate that statements concerning basic physical properties have more nontrivial models than statements concerning contingent facts of history or geography (see \cite{Heller86}).

\section{Philosophical Envoy}
At the heart of the series of works we reviewed in the previous sections is the distinction between language and what language speaks about. The importance of Tarski's achievement, in his truth definition, is that this distinction can be made without generating logical contradictions and entirely within the logical-mathematical paradigm. All the works referred to above are, in fact, exploitation of this achievement. No wonder then that the philosophical discussions around Tarski's definition are naturally transferred to our present context.

One such discussion was triggered by Hartry Field's 1972 article \cite{Field}. He argued that, in fact, Tarski did not define truth, but merely reduced some semantic terms (such as ``true'', ``denotes'', ``applies to'') to other semantic terms. In Field's opinion, a reduction of such terms to purely syntactic terms would be highly desirable as this would defend the doctrine of physicalism against its enemies. The concept of meaning, which successfully resisted all attempts of being reduced to physical concepts, remained the last safe bastion of all kinds of anti-reductionists. ``Thus in semantics -- Field writes -- physicalists agree that all legitimate semantic terms must be explicable nonsemantically -- they think in other words that there are no irreducibly semantic facts''. Facing this doctrine, ``we would have either to give up these semantic terms or else to reject physicalism''. Or perhaps, Field hopes, future developments in linguistic theory will invalidate this dichotomy.\footnote{Field's analysis of Tarski's definition of truth met with criticism of other logicians, see for instance \cite[note 7]{Etchemendy}}

Indeed, the future development of the theory of formal languages, based on categorical logic, have invalidated this dichotomy, but not quite in the sense Field expected. In all four cases discussed in the previous sections, the problem of reducing semantics to syntax has been superseded: neither of these approaches talk about reduction of semantics to syntax, or \textit{vice versa}, but rather about the adjointness between them. The difference is paramount. Adjointness establishes a specific duality between two structures, in which neither of these structures is reduced to the other, but both of them interact with each other in a mathematically strictly defined sense. First, both structures must be organized into categories (say, categories {\bf A} and {\bf B}) and then the interaction between them is described by two functors: one represents the action of category {\bf A} on category {\bf B}, the other \textit{vice versa}. The action of these adjoint functors establishes a relationship between the two categories, and does this in such a way that objects of one category are transformed into objects of the second category and are included in the network of morphisms of the latter category. The action is reciprocal: from category {\bf A} to category {\bf B}, and \textit{vice versa}.\footnote{This is by no means a definition of adjointness, but only an attempt at its intuitive description; for the definition consult any textbook on category theory.}

So far in this work we have remained exclusively within the scope of purely formal science, but -- as we indicated in the Introduction -- our analyzes may be regarded as a preparation for non-purely-formal applications. First of all, various applications related to artificial intelligence come to mind, in which the manipulation of meanings plays an important role (for some speculations in this field see \cite{AwodeyHeller}). Another field of possible applications (of enormous philosophical importance) are the wealth of issues highlighted in the fragment by Paul Davies quoted in the Introduction. He speaks about logical operations performed on organized information patterns and mathematically coded instructions, ``implying as they do semantic content''. This wording almost directly suggests the methods on which all our analyzes ware based.

Yet the purely formal aspect of the whole issue is also open to further research. For example, it would be highly desirable to investigate theories presupposing richer logical systems than those explored so far. Needless to say that with each generalization step, the indispensable mathematical tools will become more and more sophisticated. In this study, we only considered models of different theories in the the category of {\bf Sets}, that is models as functors from the syntactic category of a given theory to {\bf Sets}, but it is also fully possible to model theories in other categories. This opens up new perspectives and many possible directions for future research.

\end{document}